\documentclass[12pt]{amsart}
\usepackage{amssymb}
\usepackage{amsmath}
\usepackage{amsfonts,amsthm,latexsym,mathrsfs}
\usepackage{color}

\theoremstyle{plain}
\newtheorem{theorem}[subsection]{{\bf Theorem}}
\newtheorem{corollary}[subsection]{{\bf Corollary}}
\newtheorem{proposition}[subsection]{{\bf Proposition}}
\newtheorem{lemma}[subsection]{{\bf Lemma}}

\newtheorem{problem}[subsection]{{\bf Problem}}
\theoremstyle{remark}

\newtheorem{example}[subsection]{{\it Example}}
\numberwithin{equation}{subsection}

\def\Sym{\mathrm{Sym}}

\newcommand{\CC}{\mathbb C}

\newcommand{\PSL}{\mathrm{PSL}}

\newcommand{\PSU}{\mathrm{PSU}}

\newcommand{\Sz}{\mathrm{Sz}}

\newcommand{\im}{\mathrm{im}}
\newcommand{\PR}{\mathbb{I}}

\def \d {\mathrm d}

\def\Yn{\mathsf{\Sigma}} 

\newcommand{\Alt}{\mathrm{Alt}}
\newcommand{\Dih}{\mathrm{Dih}}

\newcommand{\ATLAS}{\textsc{Atlas\ }}

\def \Z {\mathbb Z}

\def \J{\mbox {\rm J}}
\def \Sz {\mbox {\rm Sz}}

\def \PSp {\mathrm{PSp}}
\def \PSU {\mbox {\rm PSU}}

\newcounter{ithmcount}

\newenvironment{ithm}{\begin{list}{{\rm \alph{ithmcount})}}{\usecounter{ithmcount}\labelwidth18pt
      \leftmargin18pt \topsep3pt \itemsep1pt \parsep2pt}}{\end{list}}

\makeatletter
\def\@author#1{\g@addto@macro\elsauthors{\normalsize%
    \def\baselinestretch{1}%
    \upshape\authorsep#1\unskip\textsuperscript{%
      \ifx\@fnmark\@empty\else\unskip\sep\@fnmark\let\sep=,\fi
      \ifx\@corref\@empty\else\unskip\sep\@corref\let\sep=,\fi
      }%
    \def\authorsep{\unskip,\space}%
    \global\let\@fnmark\@empty
    \global\let\@corref\@empty 
    \global\let\sep\@empty}%
    \@eadauthor={#1}
}
\makeatother

\begin{document}


\title[]
{Groups in which every non-abelian subgroup is self-normalized}
\author[C. Delizia]{Costantino Delizia}
\address{University of Salerno, Italy}
\email{cdelizia@unisa.it}
\author[U. Jezernik]{Urban Jezernik}
\address{University of Ljubljana, Slovenia}
\email{urban.jezernik@fmf.uni-lj.si}
\author[P. Moravec]{Primo\v z Moravec}
\address{University of Ljubljana, Slovenia}
\email{primoz.moravec@fmf.uni-lj.si}
\author[C. Nicotera]{Chiara Nicotera}
\address{University of Salerno, Italy}
\email{cnicoter@unisa.it}
\subjclass[2010]{}

\date{}

\begin{abstract}
We study groups having the property that every non-abelian subgroup is equal to
its normalizer. This class of groups is closely related to an open problem posed by Berkovich. 
We give a full classification of finite groups having the above
property. 
We also describe all infinite soluble groups in this class.
\end{abstract}


\maketitle

\section{Introduction}
\label{s:intro}

\noindent
In his book \cite{Ber09}, Berkovich posed the following problem:

\begin{problem}[\cite{Ber09}, Problem 9]
\label{prob:berkovich}
Study the $p$-groups $G$ in which $C_G(A)=Z(A)$ for all non-abelian $A\le G$.
\end{problem}

Classification of such $p$-groups appears to be difficult, as there seem to be many classes of finite $p$-groups enjoying the above property. 
In a recent paper \cite{Del16}, the finite $p$-groups which have maximal class or exponent $p$ and satisfy Berkovich's condition are characterized. In addition to that, the infinite supersoluble groups in which every non-abelian subgroup is self-centralized are completely classified. Some relaxations of Berkovich's problem are considered in \cite{Del13, Del15} where locally finite or infinite supersoluble groups $G$ in which every non-cyclic subgroup $A$ satisfies $C_G(A)=Z(A)$ are described.  

\smallskip

Recently, Pavel Zalesskii suggested to us another related problem:

\begin{problem}
\label{prob:zaless}
Classify finite groups $G$ in which $N_G(A)=A$ for all non-abelian $A\le G$.
\end{problem}

Let us denote by $\Yn$ the class of all groups $G$ satisfying the property stated in Problem \ref{prob:zaless}, that is, that every non-abelian subgroup of $G$ is self-normalized. Motivation for considering this problem is twofold. Firstly, every group in $\Yn$ clearly satisfies the property that every non-abelian subgroup is self-centralized. Thus the class $\Yn$ fits into the framework set by the above mentioned Berkovich's problem. Secondly, the class of groups in which every non-abelian subgroup is self-normalized can be considered as a particular case of 
the following general situation. 
Fix a property $\mathcal
P$ related to the subgroups of a given group and consider the class of all
groups $G$ in which every subgroup $H$ either has the property $\mathcal P$ or is of
bounded index in its normalizer $N_G(H)$. 
There is an abundance of literature studying restrictions these kind of conditions impose on groups.
For example, taking $\mathcal P$ to be
the property of being normal, the authors of \cite{Fer15} investigate finite
$p$-groups in this class. If we take $\mathcal P$ to be commutativity and set the upper bound of $|N_G(H):H|$ to 1, then we obtain the class $\Yn$.

\smallskip

One of the purposes of this paper is to completely characterize the finite groups in $\Yn$. This is done in Section \ref{s:finite}. We show that these
groups are either soluble or simple. Finite non-abelian simple groups in $\Yn$
are precisely the groups $\Alt(5)$ and $\PSL(2,2^{n})$, where $2^n - 1$ is a prime (see
Theorem~\ref{t:simple}). The structure of finite soluble non-nilpotent groups in
$\Yn$ is described in Theorem~\ref{t:soluble_finite}. In contrast with Berkovich's problem,
finite nilpotent non-abelian groups in
$\Yn$ are  precisely the minimal non-abelian $p$-groups, whose structure is well-known.

In Section~\ref{s:infinite} we deal with infinite groups in $\Yn$. 
We completely describe the structure of infinite soluble groups in
$\Yn$ in Theorems~\ref{t:soluble_general_case} and~\ref{t:periodic}. On the other hand, the structure of infinite soluble groups with the property that every non-abelian subgroup is self-centralized is still not completely known, the paper \cite{Del16} only deals with the supersoluble case. Finally, we
prove in Theorem \ref{t:perfect} that every infinite locally finite group in $\Yn$ is metabelian.

\section{Finite groups}
\label{s:finite}

\noindent
The following three results hold in general for groups in $\Yn$, not just for finite groups. They will be used throughout the paper without further reference.

\begin{proposition}
\label{p:basic}
The following properties hold:
\begin{enumerate}
\item The class $\Yn$ is subgroup and quotient closed.
\item If $G\in\Yn$, then every non-abelian subgroup of $G$ is self-centralized.
\item Non-abelian groups in $\Yn$ are directly indecomposable.
\end{enumerate}
\end{proposition}

\begin{proposition}
\label{p:basic2}
Let $G\in\Yn$.
\begin{enumerate}
\item Every proper subnormal subgroup of $G$ is abelian.
\item If $G$ is not perfect, then it is metabelian.
\end{enumerate}
\end{proposition}

\begin{lemma}
\label{l:soluble}
Let $G\in\Yn$ be a soluble group which is either infinite or non-nilpotent finite, and $F$ its Fitting subgroup. Then
$F$ is abelian and of prime index in $G$.
\end{lemma}
\proof
If $G$ is infinite then $F$ is abelian by \cite[Theorem 3.1]{Del16}. Otherwise, $F$ is abelian by (i) of Proposition~\ref{p:basic2}.
By  (ii) of Proposition~\ref{p:basic2} we have $G'\leq F$, so $G/F$ is an abelian group having no non-trivial proper subgroups, and hence it has prime order.
\endproof

The following easy observation immediately follows from the definition.
\begin{lemma}
\label{l:conjug}
Let $G$ be a finite group. Then $G$ belongs to $\Yn$ if and only if the following holds: for every non-abelian subgroup $H$ of $G$, the number of conjugates of $H$ in $G$ is equal to $|G:H|$.
\end{lemma}

The next lemma gives a well-known characterization of minimal non-abelian finite $p$-groups, see \cite[Exercise 8a, p.\ 29]{Ber09}.

\begin{lemma}
\label{l:minnonab}
Let $G$ be a finite $p$-group. The following are equivalent:
\begin{ithm}
\item $G$ is minimal non-abelian.
\item $\d(G)=2$ and $|G'|=p$.
\item $\d(G)=2$ and $Z(G)=\Phi(G)$.
\end{ithm}
\end{lemma}

\begin{proposition}
\label{p:pgroups}
Let $G$ be a nilpotent group. Then $G\in\Yn$ if and only if $G$ is either abelian or a finite minimal non-abelian $p$-group for some prime $p$. 
\end{proposition}

\proof
Let $G$ be a nilpotent group in $\Yn$, and suppose that $G$ is not abelian.
By \cite[Theorem 3.1]{Del16}, $G$ has to be finite, and Proposition~\ref{p:basic} implies that $G$ is a $p$-group. Let $H$ be a proper subgroup of $G$. Then $H<N_G(H)$ by \cite[5.2.4]{Rob96}, hence $H$ is abelian. This proves that $G$ is minimal non-abelian. The converse is clear.
\endproof

Finite soluble groups in which all Sylow subgroups are abelian are called {\em $A$-groups}, cf. \cite[Seite 751]{Hup67}. Next we show that finite soluble non-nilpotent groups in $\Yn$ are $A$-groups. 

\begin{lemma}
\label{l:soluble_finite}
Let $G \in \Yn$ be a soluble non-nilpotent group and let $p = |G : F|$. Then the
Sylow $p$-subgroups of $G$ are cyclic and, for primes $r \neq p$,  the Sylow
$r$-subgroups are abelian.
\end{lemma}
\proof
If $r \neq p$, then the Sylow $r$-subgroups of $G$ are contained in $F$ and so
are abelian by \ref{l:soluble}.
Now let $P$ be a Sylow $p$-subgroup of $G$ and assume that $P$ is not cyclic.
Since $C_G(F) = F$, there exists a prime $r$ and a Sylow $r$-subgroup $R$ of $G$
such that $P$ and $R$ do not commute. Since $R$ is normal in $G$, $X=P R$ is a
subgroup of $G$. Since $P$ is not cyclic, there are maximal subgroups $P_1$ and
$P_2$ of $P$ with $P_1 \neq P_2$. It follows that $X$ normalizes both $P_1R$ and
$P_2R$. As $X > P_1R$ and $X > P_2R$, we conclude that $P_1R$ and $P_2R$ are
abelian. But then $P = P_1P_2$ commutes with $R$, a contradiction. Hence $P$ is
cyclic.
\endproof

Suppose that an element $x$ acts on an abelian group $A$. Consider the induced
homomorphism $\partial_x \colon A \to A$, $\partial_x(a) = a^{1-x}$. We will
describe groups belonging to the class $\Yn$ based on the following property of
$\partial_x$:
\begin{equation}
\label{eq:I}
\forall B \leq A\colon \; ( \partial_x (B) \leq B \; \Longrightarrow \; \partial_x (B) = B).
\tag{$\PR$}
\end{equation}
Note that $x$ acts fixed point freely on $A$ if and only if $\partial_x$ is
injective. The property $\PR$ implies injectivity of $\partial_x$, since
taking $B = \ker \partial_x$ immediately gives $B = 1$. Furthermore, the
property $\PR$ implies that $\partial_x$ is an epimorphism by taking $B = A$.
Therefore having property $\PR$ implies that $\partial_x$ is an isomorphism.

The following proposition shows how property $\PR$ is related to $\Yn$.

\begin{proposition}
\label{p:property}
Let $G = \langle x \rangle \ltimes A$ with $x^p$ acting trivially on $A$ for
some prime $p$ and $x$ acting fixed point freely on $A$. 
Then $G \in \Yn$ if and only if $\partial_x$ has property $\PR$.
\end{proposition}
\proof
Assume first that $G \in \Yn$. If $\partial_x$ does not have property $\PR$,
then there exists a subgroup $B \leq A$ such that $\partial_x(B) \subsetneq B$.
Consider the subgroup $H = \langle x \rangle \ltimes \partial_x(B)$ of $G$. Take
an element $b \in B \backslash \partial_x(B)$ and observe that $x^b = x
\partial_x(b) \in H$. Therefore $b$ normalizes $H$ and does not belong to $H$.
This implies that $H$ is abelian, and so $\partial_x(\partial_x(B)) = 1$.
By injectivity of $\partial_x$, it follows that $B = 1$, a contradiction with
$\partial_x(B) \subsetneq B$. Therefore $\partial_x$ has property $\PR$.

Conversely, assume now that $\partial_x$ has property $\PR$. To prove that $G$
belongs to $\Yn$, take a non-abelian subgroup $H$ of $G$. Note that $H$ must
contain an element of the form $x a$ for some $a \in A$. Then since $\partial_x$
is surjective, we have $a = \partial_x(b)$ and so $xa = x^b$. After possibly
replacing $H$ by $H^{b^{-1}}$, it suffices to consider the case when $b = 1$,
and therefore $x \in H$. We can thus write $H = \langle x \rangle  \ltimes B$
for some $B = H \cap A \leq A$. Let us show that $H$ is self-normalizing in $G$.
To this end, take an element $x^j c \in N_G(H)$. Then  $x^{x^j c} = x^c = x
c^{-x} c$, and so we must have $\partial_x(c) \in B$. Conversely, any element
$x^j c \in A$ with the property that $\partial_x(c) \in B$ normalizes $H$, since
for any $b \in B$ we also have $b^{x^j c} = (b^{x^j})^c = b^{x^j} \in B$. Thus
$N_G(H) = \langle x \rangle \ltimes  \partial_x^{-1}(B)$. By
property $\PR$, we have $\partial_x^{-1}(B) =
B$, which implies $N_G(H) = \langle x \rangle \ltimes B = H$, as required.
\endproof

\begin{example} \label{ex:group_ring}
Let $\zeta_p$ be a primitive complex $p$-th root of unity for some prime $p$.
Then $\zeta_p$ acts by multiplication on the abelian group $\CC$ and we can form
$G = \langle \zeta_p \rangle \ltimes \CC$. The group $G$ has a subgroup $H =
\langle \zeta_p \rangle \ltimes \Z[\zeta_p]$. Now, note that
$\zeta_p$ acts fixed point freely on $\CC$, and multiplication by $1- \zeta_p$
is invertible on $\CC$. Therefore $\partial_{\zeta_p}$ is an isomorphism of
$\CC$. However, multiplication by $1 - \zeta_p$ maps $\Z[\zeta_p]$ into its
augmentation $\{ \sum_i \lambda_i \zeta_p^i \mid \sum_i = 0 \}$. In particular, the restriction of
$\partial_{\zeta_p}$ on $\Z[\zeta_p]$ is not surjective. Therefore the isomorphism
$\partial_{\zeta_p}$ does not have property $\PR$. In particular, neither $G$
nor $H$ belong to the class $\Yn$.
\end{example}

\begin{example}
Take $C_p = \langle x \rangle$ for some prime $p$ and consider a
$\Z[C_p]$-module $A$. Assume that $\partial_x$ is an isomorphism
of $A$. To verify whether or not $\partial_x$ has property $\PR$, it suffices
to show that the restriction of $\partial_x$ on every cyclic submodule of $A$
is surjective. To this end, suppose that $B$ is a cyclic $\Z[C_p]$-module with
$\partial_x$ having trivial kernel on $B$.
Therefore $B$ is isomorphic to a quotient of the ring $\Z[C_p]$ by some ideal
$J$. Denote
$D = x - 1$ and $N = x^{p-1} + x^{p-2} + \cdots + 1$ as elements in $\Z[C_p]$.
Observe that injectivity of $\partial_x$ is equivalent to saying that
whenever $D \cdot z \in J$ for some $z \in B$, it follows that $z \in J$. 
Now, we have that $D \cdot N = 0$, and so it follows that $N \in J$. 
Therefore $J$ is the preimage of an ideal in the ring $\Z[C_p]/N\Z[C_p] \cong
\Z[\zeta_p]$, where $\zeta_p$ is a complex primitive $p$-th root of unity.
Note that multiplication by $D$
is surjective on $\Z[C_p]/J$ if and only if we have $\im D + J = \Z[C_p]$.

Consider two special cases. First let $J = N \Z[C_p]$. This corresponds to the
module $\Z[\zeta_p]$ from Example \ref{ex:group_ring}.  Since $\im N = \ker D$
(see \cite[Lemma 9.26]{Rot09}), we have that $\im D + J = \im D + \im N$.
Dividing the polynomial $N$ by $D$ in $\Z[\langle x \rangle]$, we get the
remainder $p \in \Z$. Whence $\im D + \im N$ contains  $p \Z[C_p]$. On the other
hand, $\im D + \im N$ is not the whole of $\Z[C_p]$, since $\partial_{\zeta_p}$
is not surjective on $\Z[\zeta_p]$. Consider now the case when $J$ is the ideal
generated by $N$ and a prime $q$ distinct from $p$. Thus $\im D + J$ contains
$p\Z[C_p]$ and $q$. It follows that $\im D + J = \Z[C_p]$, and $\partial_x$ is
surjective in this case. Moreover, the map $\partial_x$ will be surjective on
all cyclic submodules of $\Z[C_p]/J$, as the ideals corresponding to these
submodules all contain $J$. 
Therefore
the group $\langle x \rangle \ltimes \Z_q[\zeta_p]$ belongs to $\Yn$.
\end{example}

\begin{lemma} \label{l:p_divisible}
Let $x$ be an automorphism of order $p$ of an abelian group $A$. If $\partial_x$
is surjective, then $A = pA$, i.e., $A$ is $p$-divisible.
\end{lemma}
\proof
Consider $A$ as a $\Z[\langle x \rangle]$-module. In this sense, the operator
$\partial_x$ corresponds to the element $1-x \in \Z[\langle x \rangle]$. We
have $(1-x)^p \equiv 0$ modulo $p\Z[\langle x \rangle]$, and so the image of
$(\partial_x)^p \colon A \to A$ is a subgroup of $p \Z[\langle x \rangle] A = pA$.
As $\partial_x$ is assumed to be surjective, it follows that $A = pA$.
\endproof

\begin{corollary}
\label{c:free_abelian}
Let $G = \langle x \rangle \ltimes A$ with $x^p$ acting trivially on $A$ for
some prime $p$. Assume that $A$ is free
abelian of finite rank. Then $G$ does not belong to $\Yn$.
\end{corollary}
\proof
By Lemma \ref{l:p_divisible}, the map $\partial_x$ is not surjective,
and so $\partial_x$ does not have property $\PR$. 
It follows from Proposition~\ref{p:property} that $G$ does
not belong to $\Yn$.
\endproof

\begin{theorem}
\label{t:soluble_finite}
Let $G$ be a finite soluble non-nilpotent group. Then $G \in \Yn$ if and only
if $G$ splits as $G = \langle x \rangle \ltimes A$, where $\langle x \rangle$ is
a $p$-group for some prime $p$, $A$ is an abelian $p'$-group, $x^p$ is central and
$x$ acts fixed point freely on $A$. 
\end{theorem}
\proof

Assume first that $G \in \Yn$. By Lemma~\ref{l:soluble_finite}, all Sylow subgroups of $G$ are abelian.
It follows that $G' \cap Z(G) = 1$ by \cite[10.1.7]{Rob96}, and $G$ splits as $G = \langle x
\rangle \ltimes G'$ with $x^p$ in the Fitting subgroup of $G$ for some prime
$p$. Whence $\langle x^p \rangle \leq Z(G)$. 
 If an element $x a\in G$ is
central, then $a^x = a$ and so $a$ must be central. As $x$ is not central, we
must have $Z(G) = \langle x^p \rangle$ and $C_G(x) = \langle x \rangle$.
Observe that as $G$ belongs to $\Yn$, the map $\partial_x$ is surjective on $G'$,
and so by Lemma \ref{l:p_divisible} the group $G'$ must be of $p'$-order.
Now, if $\langle x \rangle$ is not of prime power order, then it splits as a
product $\langle x \rangle = A_p \times A_{p'}$ with $A_p$ a $p$-group and $A_{p'}$
a $p'$-group. Then $A_p \ltimes G'$ is a non-abelian proper normal subgroup of
$G$, a contradiction. Whence $\langle x \rangle$ is of $p$-power order.
Note that $x$ acts fixed
point freely on $G'$ since
$C_G(x) \cap G' = 1$. Thus $G/Z(G)$ is a Frobenius group with complement
of order $p$.

Conversely, take $G = \langle x \rangle \ltimes A$ with the stated properties.
Therefore $x^p$ acts trivially on $A$ and $\partial_x$ is an injective
endomorphism of $A$. As $A$ is finite,  it immediately follows that $\partial_x$
is surjective and that it satisfies property $\PR$. It now follows from Lemma
\ref{p:property} that $G$ belongs to $\Yn$.
\endproof

Notice that in Theorem~\ref{t:soluble_finite} we have that $\langle x\rangle$ is the Sylow $p$-subgroup of $G$ and $F=\langle x^p\rangle\ltimes A$, so $A=F/\langle x^p\rangle$.

\begin{theorem}
\label{t:solvsimp}
Let $G$ be a finite group in $\Yn$. Then $G$ is either soluble or simple.
\end{theorem}

\proof 
By induction on the order of $G$. Suppose that $G$ is not simple, and let $A$ be a maximal normal subgroup of $G$. Then $A$ is abelian since $G$ is in $\Yn$, hence $C_G(A)$ contains $A$ and is normal in $G$. It follows that either $C_G(A)$ is abelian or $C_G(A)=G$.

Assume first that $C_G(A)$ is abelian. Thus $A=C_G(A)$ by the maximality. Let $P/A$ be a Sylow subgroup of $G/A$, and choose $x\in P\setminus A$. Then $A\langle x\rangle$ is not abelian and subnormal in $P$, so $A\langle x\rangle=P$ by (i) of Proposition~\ref{p:basic2}. This shows that every Sylow subgroup of $G/A$ is cyclic, hence $G/A$ is soluble. Therefore $G$ is soluble.

Hence we can assume that $C_G(A)=G$, and so $A\leq Z(G)$. Moreover, $A$ is contained in every maximal subgroup $M$ of $G$. Namely, if this is not the case, then $G=MA$ implies that $G'=M'<G$, so $G'$ is abelian, hence $G$ is soluble and we are done. Therefore every maximal subgroup of $G$ is non-simple, so it is soluble by the induction hypothesis. Suppose that $G$ has a maximal subgroup which is nilpotent. Then Proposition~\ref{p:pgroups} implies that its Sylow $2$-subgroup has nilpotency class $\leq 2$, therefore it follows from \cite{Jan61} that $G$ is soluble. Hence we can assume that every maximal subgroup of $G$ is not nilpotent. In this case, every Sylow subgroup of $G$ is abelian by Lemma~\ref{l:soluble_finite}. Then we have $G'\cap Z(G)=1$, hence $G'<G$. Therefore $G'$ is abelian, $G$ is soluble and the proof is complete.
\endproof 

Theorem~\ref{t:solvsimp}, Proposition~\ref{p:pgroups} and Theorem~\ref{t:soluble_finite} show that, in order to obtain a full classification of all finite groups in the class $\Yn$, it only remains to describe the finite simple groups in $\Yn$. At first we need a couple of auxiliary results.

\begin{lemma}
\label{l:dihedral}
Let $n>2$. 
The dihedral group $\Dih(n)$ of order $2n$ belongs to $\Yn$ if and only either $n=4$ or $n$ is odd.
\end{lemma}

\proof
Denote  $G=\Dih(n)=\langle x,y\mid x^n=y^2=1,x^y=x^{-1}\rangle$, where $n\ne4$.
Suppose that $G\in\Yn$. Then $n$ is not a power of 2 by Proposition~\ref{p:pgroups}. Let $p$ be an odd prime dividing $n$, and assume that $n$ is even. 
Denote $H=\langle y,x^{n/p}\rangle$.
Then
$x^{n/2}\in Z(G)\setminus H$, hence $H$ is not self-normalized. Therefore $n$ is odd.

Conversely, clearly $\Dih(4)$ belongs to $\Yn$ since it is minimal non-abelian. Suppose now $n$ is odd. Let $H$ be a non-abelian subgroup of $\Dih(n)$ of index $m$. Then $H$ is conjugate to $K=\langle x^m,y\rangle$. Take $z=x^iy^j\in N_G(K)$, $0\le i< n$, $0\le j\le 1$. Then
$y^{x^iy^j}=x^{2(-1)^{j+1}i}y\in K$ if and only $m$ divides $i$, that is, $z\in K$. This shows that $\Dih(n)\in\Yn$.
\endproof

\begin{lemma}
\label{l:psl}
If $q\neq 3,5$ is an odd prime power, then $\PSL(2,q)$ does not belong to $\Yn$.
\end{lemma}

\proof
Let $G=\PSL(2,q)$. Since $q$ is odd, it follows from \cite{Dic01} that $G$ contains a subgroup $H$ isomorphic to $\Dih ((q-1)/2)$, and a subgroup $K$ isomorphic to $\Dih ((q+1)/2)$. If $q\equiv 1\mod 4$, then $H$ is not in $\Yn$, unless $q=5$, whereas if $q\equiv 3\mod 4$, then $K\notin\Yn$, unless $q=3$ or $q=7$. Notice that $\PSL(2,7)\notin\Yn$ since it has a subgroup isomorphic to $\Sym(4)$ (see \cite[Theorem 8.27]{Hup67}), and therefore  a non-abelian subgroup which is isomorphic to $\Alt(4)$ and is not self-normalizing.
\endproof

\begin{theorem}
\label{t:simple}
A finite non-abelian simple group $G$ belongs to $\Yn$ if and only if it is isomorphic to $\Alt (5)$ or $\PSL(2,2^n)$, where $2^n - 1$ is a prime.
\end{theorem}

\proof
Let $G\in\Yn$ be finite non-abelian simple. Let $P_p$ be a Sylow $p$-subgroup of $G$. It follows from \cite[Theorem 1.1]{Gur03} that if $p>3$, then $P_p$ is abelian. For $p=3$, the same result implies that either $P_3$ is abelian or $G\cong \PSL (2, 3^{3^a})$, where $a\ge 1$. The latter cannot happen by Lemma~\ref{l:psl}. Hence we conclude that $P_3$ needs to be abelian.
If $P_2$ is also abelian, then all Sylow subgroups of $G$ are abelian, and it follows from \cite{Bro71} that $G$ belongs to one of the following groups:
$\J_{1}$, or $\PSL(2,q)$, where $q>3$ and $q\equiv 0,3,5 \mod 8$. Note that the latter condition can be reduced to $q\equiv 0\mod 8$ or $q=5$ by Lemma~\ref{l:psl}.
If $P_2$ is non-abelian, then it is minimal non-abelian by Proposition~\ref{p:pgroups}, and hence $P_2$ is nilpotent of class two. By \cite{Gil75}, $G$ is isomorphic to one of the following groups: $\PSL(2,q)$, where $q\equiv 7,9 \mod 16$, $\Alt(7)$, $\Sz(2^n)$, $\PSU (3,2^n)$, $\PSL (3,2^n)$ or $\PSp (4,2^n)$, where $n\ge 2$. The first family can be ruled out by Lemma~\ref{l:psl}.

It suffices to see which of the above listed groups belong to $\Yn$. It follows from \ATLAS that the Janko group $\J_1$ has a subgroup isomorphic to $\Dih(3)\times \Dih(5)$, hence it is not in $\Yn$ by Proposition~\ref{p:basic}. Also, note that $\Alt(7)$ has a subgroup isomorphic to $\Sym(4)$, therefore $\Alt(7)\notin\Yn$.
If $G=\Sz (2^n)$ and $P_2$ is its Sylow 2-subgroup, then $|P_2'|=2^n$ by \cite{Suz62}, hence the Suzuki groups do not belong to $\Yn$. Similarly, if $G$ is $\PSU (3,2^n)$ or $\PSL (3,2^n)$, then the derived subgroup of a Sylow 2-subgroup of $G$ has order $2^n$, whereas if $G=\PSp (4,2^n)$, then $|P_2'|=2^{2n}$. This shows that neither of these groups belongs to $\Yn$.

We are left with the groups $G=\PSL(2,q)$, where $q=5$ or $q\equiv 0 \mod 8$. The subgroup structure of $G$ is described in \cite{Dic01}. 
It is straightforward to verify
that $\PSL (2,5)=\Alt(5)\in \Yn$. 
Consider now $q=2^n$ with $n\ge 3$. 
Suppose first $q - 1$ is not a prime. Let $d$ be a proper divisor of $q - 1$.
Then it follows from \cite[section 250]{Dic01} that $G$ has a single conjugacy
class of size $q + 1$ of a subgroup $H$ isomorphic to $C_2^n\rtimes C_d$.
Therefore $|G:N_G(H)|=q+1$, and $|G:H|=(q^2-1)/d$. As $d < q - 1$, it follows
that $H\neq N_G(H)$, therefore $G\notin\Yn$.
On the other hand, let $q - 1$ now be a prime. Going through the list of
subgroups of $G$ given in \cite{Dic01}, along with the given data on the number
of conjugates of these subgroups, we see that apart from the subgroups in
section 250, one has that for every non-abelian subgroup $H$ of $G$, the number
of conjugates of $H$ is equal to $|G:H|$. As for the remaining subgroups, note
that they must be of order $2^m d$ for some integer $m$ and some divisor $d$ of
$q - 1$. There are only two such options, one corresponding to abelian groups of
order $2^m$ and the other corresponding to subgroups of order $q(q-1)$. Each of
the latter belongs to a system of $(q^2 - 1)2^{n-m} / (2^k - 1)$ conjugate
groups for some $k$ dividing $n$. Note that since $q-1$ is a prime,  $n$ must
also be a prime. When $k = 1$,
the group under consideration is abelian; therefore $k = n$. We thus have that  the number of conjugates of
each of these subgroups is equal to their index in $G$. By Lemma~\ref{l:conjug},
this shows that if $2^n - 1$ is a prime, then $\PSL(2,2^n)$ belongs to $\Yn$.
\endproof


\section{Infinite groups}
\label{s:infinite}

Let $G\in\Yn$ be an infinite finitely generated soluble group, and let $F$
denote the Fitting subgroup of $G$. Then $F$ is polycyclic and $G=\langle
x\rangle F$ for every element $x\in G\setminus F$, by Lemma~\ref{l:soluble}. We
will denote by $h(F)$ the Hirsch length of $F$. In what follows, the set of all
periodic elements of a group $G$ will be denoted by $T(G)$. For a prime $p$,
let $T_p(G)$ be the set of elements in $G$ of $p$-power order, and let $T_{p'}(G)$
be the set of elements in $G$ of order coprime to $p$.

\begin{lemma}
\label{l:hirsch}
Let $G\in\Yn$ be an infinite finitely generated soluble group, and suppose $h(F)=1$. Then $G$ is abelian.
\end{lemma}
\proof Assume not. Since $h(F)=1$ the group $G$ is infinite cyclic-by-finite. It easily follows that there exists a finite normal (and hence abelian) subgroup $N$ such that $G/N$ is either infinite cyclic or infinite dihedral. As the latter group is not in $\Yn$, we can write $G=\langle x\rangle\ltimes N$ where $x$ aperiodic. Since $N\leq F$ and $G$ is not abelian, we conclude that $x\notin F$. By Lemma~\ref{l:soluble} there exists a prime number $p$ such that $x^p\in F$. Then $x^p\in Z(G)$.  Then $G/\langle x^p\rangle$ is finite. Moreover, it is not nilpotent by \cite[Theorem 3.1]{Del16}. The Fitting subgroup of  $\frac G{\langle x^p\rangle}$ is  $\frac F{\langle x^p\rangle}=\frac{\langle x^p\rangle\times T(F)}{\langle x^p\rangle}$, so it equals $T(F)$ since $h(F)=1$. Now by Theorem~\ref{t:soluble_finite} it follows that $\frac G{\langle x^p\rangle}=\frac{\langle x \rangle}{\langle x^p \rangle} \ltimes T(F)$. Hence $G=\langle x\rangle\ltimes T(F)$. Let $q\ne p$ be any prime number which does not divide the order of $T(F)$. Thus $\frac{\langle x^q \rangle}{\langle x^{pq} \rangle} \ltimes T(F)$ is a normal proper subgroup of the factor group  $G/\langle x^{pq}\rangle$. Note that $\frac{\langle x^q \rangle}{\langle x^{pq} \rangle} \ltimes T(F)$ is not abelian since $x\notin Z(G)$. Therefore $G/\langle x^{pq}\rangle\notin\Yn$, a contradiction.
\endproof

\begin{theorem}
\label{t:soluble_general_case}
Let $G\in\Yn$ be an infinite finitely generated soluble group. Then $G$ is
abelian.
\end{theorem}
\proof
We have $G = \langle x \rangle F$ with  $x^p \in F \cap Z(G)$. Consider first
the case when $x^p = 1$. Thus $G = \langle x \rangle \ltimes F$.

Let us show that $x$ acts fixed point freely on $F$. To this end, let $f \in F$
be an element with $f^x = f$. For any positive integer $k$, the quotient groups
$G/F^k$ are finite and belong to $\Yn$. This shows that $f^x F^k = f F^k$, and thus
Theorem~\ref{t:soluble_finite} gives that $f F^k$ is trivial in $G/F^k$. This
means that $f \in \bigcap_{k} F^k = 1$. Therefore $x$ acts fixed point freely on
$F$ and $\partial_x$ is injective.

Since $G \in \Yn$, the group $F$ can not be free abelian by  Corollary
\ref{c:free_abelian}. Thus the torsion subgroup $T(F)$ is not trivial. The
factor group $G/T(F)$ belongs to $\Yn$, and so the action of $x$ on this group
is trivial by Corollary~\ref{c:free_abelian}. Therefore $\partial_x(F)$ is trivial
in $G/T(F)$. So the image of the map $\partial_x$ is contained in the finite
group $T(F)$. This is impossible since $\partial_x$ is injective and its
domain $F$ is infinite.

Lastly, consider the case when $x^p$ is not trivial in $G$. Let us look at the group
$G/\langle x^p \rangle$. If $h(G) = 1$, then $G$ is abelian by Lemma
\ref{l:hirsch}. Thus $h(G) \geq 2$ and so $G/\langle x^p \rangle$ is an infinite
finitely generated soluble group with an element of order $p$ outside its
Fitting subgroup. Therefore $G/ \langle x^p \rangle$ must be abelian by the above
argument. As $\langle x^p \rangle$ is contained in $Z(G)$, it follows that $G$
is nilpotent, and so $G$ must be abelian. The proof is now complete.
\endproof

Let $G\in \Yn$ be a soluble group. It follows easily from
Theorem~\ref{t:soluble_general_case} that every aperiodic element of $G$ is
central. As a consequence, every non-periodic soluble group in $\Yn$ is abelian.
Therefore the following result completes the description of all infinite soluble
groups in $\Yn$.

\begin{theorem}
\label{t:periodic}
Let $G$ be an infinite non-abelian soluble periodic group. Then $G$ belongs
to $\Yn$ if and only if $G$ splits as $\langle x \rangle \ltimes A$ with
$\langle x \rangle$ a $p$-group for some prime $p$, $A$ is a $p'$-group,
$x^p$ acts trivially on $A$ and $\partial_x$ has property $\PR$.
\end{theorem}
\proof
It follows from Proposition~\ref{p:property} that a group $G$ with the above
decomposition belongs to $\Yn$. Therefore we are only concerned with proving
the converse.

Assume that $G \in \Yn$ is an infinite non-abelian soluble periodic group.
We have that $G = \langle x \rangle F$ with $x^p \in F \cap Z(G)$ for some prime
$p$. Note that since $x$ is of finite order, say, $p^k \beta$ for some $\beta$
coprime to $p$, we may replace $x$ by $x^\beta$ and assume from now on that
$\langle x \rangle \subseteq T_p(F)$.

Let us first prove that $T_p(F) = \langle x^p \rangle$. To this end, it suffices
to consider the factor group $G / \langle x^p \rangle$, and therefore we can
assume that $x^p = 1$. Thus $G = \langle x \rangle \ltimes F$, and so $G =
(\langle x \rangle \ltimes T_p(F)) \ltimes T_{p'}(F)$. If the group $\langle x
\rangle \ltimes T_p(F)$ is not cyclic, then there is an element $z \in T_p(F)$
that commutes with $x$. It follows that the group $\langle x
\rangle \ltimes T_p(F)$ contains the subgroup $\langle x
\rangle \ltimes \langle z \rangle \cong C_p \times C_p$. Now $G$ contains the
subgroup $\langle x \rangle \ltimes T_p(F)$ that is normalized by the element
$z$. This is a contradiction with $G \in \Yn$. Hence the group $\langle x
\rangle \ltimes T_p(F)$ is cyclic. This is possible only if $T_p(F)$ is trivial,
as claimed.

We now have a splitting $G = \langle x \rangle \ltimes A$ with $A = T_{p'}(F)$,
$x$ acts non-trivially on $A$ and $x^p$ is central. Let us now show that $x$ acts
fixed point freely on $A$. It will then follow from Proposition~\ref{p:property}
that $\partial_x$ has property $\PR$. To see this, assume that $z \in A$ is a
fixed point of $x$. Thus $z \in Z(G)$. In particular, as $G \in \Yn$, we have
that $z$ must be contained in every non-abelian subgroup of $G$. Now, as $x$ acts
non-trivially on $A$, there is an element $b \in A$ with $b^x \neq b$. Set $B =
\langle b, b^x, \dots, b^{x^{p-1}} \rangle$, this is an $x$-invariant finite
subgroup of $G$. Therefore $G$ possesses the finite non-abelian subgroup $\langle
x \rangle \ltimes B$. By Theorem~\ref{t:soluble_finite}, we have that $x$ acts
fixed point freely on $B$. On the other hand, we must have that $z \in \langle x
\rangle \ltimes B$, and so $z \in B$. This implies that $z$ is trivial, as
required. The proof is complete.
\endproof

Let $G$ be an infinite group in $\Yn$, and suppose that $G$ is not soluble. Then $G$ is perfect by (ii) of Proposition~\ref{p:basic2}. Our last result gives information on the structure of such a group $G$ provided that it is locally finite.

\begin{theorem}
\label{t:perfect}
Let $G$ be an infinite locally finite group in $\Yn$. Then $G$ is metabelian.
\end{theorem}

\proof
Let $G\in\Yn$ be locally finite, and suppose that $G$ is not metabelian. Then
$G$ contains a finite insoluble subgroup, say $H_0$.  It follows from
Theorem~\ref{t:solvsimp} and Theorem~\ref{t:simple} that $H_0$ is isomorphic to
either $\Alt(5)$ or some $\PSL(2,2^n)$ with $n$ a prime. Pick an element $x_1
\in G$ that does not belong to $H_0$, and set $H_1 = \langle x_1, H_0 \rangle$.
We now have that the group $H_1$ is isomorphic to some $\PSL(2,2^m)$ with $m$ a
prime. Finally let $x_2 \in G$ be an element not in $H_1$, and set $H_2 =
\langle x_2, H_1 \rangle$. The group $H_2$ is isomorphic to some $\PSL(2,2^k)$
with $k$ a prime. Now, as $H_2 = \PSL(2,2^k)$ properly contains $H_1 =
\PSL(2,2^m)$, we must have $m \mid k$, which is impossible since $m$ and $k$ are
distinct primes.
\endproof


\end{document}